\documentclass[12pt,vatola]{article}

\usepackage{graphicx}

\def\c{\centerline}

\def\re#1{\par\hangindent\parindent\indent\llap{#1\enspace}\ignorespaces}

\def\no{\noindent}

\begin{document}

\c{\Large\bf Smarandache Multi-Space Theory(IV) }\vskip 5mm

\hskip 50mm {\it -Applications to theoretical physics}\vskip 10mm

\c{Linfan Mao}\vskip 3mm

\c{\small Academy of Mathematics and System Sciences}

\c{\small Chinese Academy of Sciences, Beijing 100080}

\c{\small maolinfan@163.com}

\vskip 10mm

\begin{minipage}{130mm}

\no{\bf Abstract.} {\small A Smarandache multi-space is a union of
$n$ different spaces equipped with some different structures for
an integer $n\geq 2$, which can be both used for discrete or
connected spaces, particularly for geometries and spacetimes in
theoretical physics. This monograph concentrates on characterizing
various multi-spaces including three parts altogether. The first
part is on {\it algebraic multi-spaces with structures}, such as
those of multi-groups, multi-rings, multi-vector spaces,
multi-metric spaces, multi-operation systems and multi-manifolds,
also multi-voltage graphs, multi-embedding of a graph in an
$n$-manifold,$\cdots$, etc.. The second discusses {\it Smarandache
geometries}, including those of map geometries, planar map
geometries and pseudo-plane geometries, in which the {\it Finsler
geometry}, particularly the {\it Riemann geometry} appears as a
special case of these Smarandache geometries. The third part of
this book considers the {\it applications of multi-spaces to
theoretical physics}, including the relativity theory, the
M-theory and the cosmology. Multi-space models for $p$-branes and
cosmos are constructed and some questions in cosmology are
clarified by multi-spaces. The first two parts are relative
independence for reading and in each part open problems are
included for further research of interested readers.}

\vskip 3mm \no{\bf Key words:} {\small multi-space, relativity
theory, M-theory, $p$-brane, multi-space model of cosmos.}

 \vskip 3mm \no{{\bf
Classification:} AMS(2000) 03C05,05C15,51D20,51H20,51P05,83C05,
83E50}
\end{minipage}

\newpage

\no{\large\bf Contents}\vskip 8mm

\no $6.$ Applications to theoretical physics\dotfill 3\vskip 3mm

\no $\S 6.1$ \ Pseudo-Faces of Spaces\dotfill 3\vskip 2mm

\no $\S 6.2$ \ Relativity Theory\dotfill 7\vskip 2mm

\no $\S 6.3$ \ A Multi-Space Model for Cosmoses\dotfill 11\vskip
2mm

\no $6.3.1$ What is M-theory\dotfill 11

\no $6.3.2$ A pseudo-face model for $p$-branes \dotfill 14

\no $6.3.2$ A multi-space model of cosmos \dotfill 18\vskip 3mm

\no  References \dotfill 21

\newpage

\no{\large\bf $6.$ Applications to theoretical physics}

\vskip 35mm

\no Whether are there finite, or infinite cosmoses? Is there just
one? What is the dimension of our cosmos? Those simpler but more
puzzling problems have confused the eyes of human beings thousands
years and one does not know the answer even until today. The
dimension of the cosmos in the eyes of the ancient Greeks is $3$,
but Einstein's is $4$. In recent decades, $10$ or $11$ is the
dimension of our cosmos in superstring theory or M-theory. All
these assumptions acknowledge that there is just one cosmos. Which
one is the correct and whether can human beings realize the cosmos
or cosmoses? By applying results gotten in Chapters $3$-$5$, we
tentatively answer those problems and explain the Einstein's or
Hawking's model for cosmos in this chapter.

\vskip 8mm

\no{\bf \S $6.1$ \ Pseudo-Faces of Spaces}

\vskip 3mm

\no Throughout this chapter, ${\bf R}^n$ denotes an Euclid space
of dimensional $n$. In this section, we consider a problem related
to how to represent an Euclid space in another. First, we
introduce the conception of pseudo-faces of Euclid spaces in the
following.

\vskip 4mm

\no{\bf Definition $6.1.1$} \ {\it Let ${\bf R}^m$ and $({\bf
R}^n,\omega)$ be an Euclid space and a pseudo-metric space. If
there is a continuous mapping $p: {\bf R}^m\rightarrow ({\bf
R}^n,\omega)$, then the pseudo-metric space $({\bf
R}^n,\omega(p({\bf R}^m)))$ is called a pseudo-face of ${\bf R}^m$
in $({\bf R}^n,\omega)$.}

\vskip 3mm

Notice that these pseudo-faces of ${\bf R}^3$ in ${\bf R}^2$ have
been considered in Chapter $5$. For the existence of a pseudo-face
of an Euclid space ${\bf R}^m$ in ${\bf R}^n$, we have a result as
in the following.

\vskip 4mm

\no{\bf Theorem $6.1.1$} \ {\it Let  ${\bf R}^m$ and $({\bf
R}^n,\omega)$ be an Euclid space and a pseudo-metric space. Then
there exists a pseudo-face of ${\bf R}^m$ in $({\bf R}^n,\omega)$
if and only if for any number $\epsilon >0$, there exists a number
$\delta >0$ such that for $\forall
\overline{u},\overline{v}\in{\bf R}^m$ with
$\|\overline{u}-\overline{v}\|<\delta$,}

$$\|\omega(p(\overline{u}))-\omega(p(\overline{v}))\|<\epsilon,$$

\no{\it where $\|\overline{u}\|$ denotes the norm of a vector
$\overline{u}$ in the Euclid space.}

\vskip 3mm

{\it Proof} \ We only need to prove that there exists a continuous
mapping $p: {\bf R}^m\rightarrow ({\bf R}^n,\omega)$ if and only
if all of these conditions in this theorem hold. By the definition
of a pseudo-space $({\bf R}^n,\omega)$, since $\omega$ is
continuous, we know that for any number $\epsilon >0$,
$\|\omega(\overline{x})-\omega(\overline{y})\| <\epsilon$ for
$\forall\overline{x},\overline{y}\in{\bf R}^n$ if and only if
there exists a number $\delta_1
>0$ such that $\|\overline{x}-\overline{y})\|<\delta_1$.

By definition, a mapping $q:{\bf R}^m\rightarrow{\bf R}^n$ is
continuous between Euclid spaces if and only if for any number
$\delta_1 >0$, there exists a number $\delta_2 >0$ such that
$\|q(\overline{x})-q(\overline{y})\|<\delta_1$ for $\forall
\overline{u},\overline{v}\in {\bf R}^m$ with
$\|\overline{u}-\overline{v})\|<\delta_2$.

Combining these assertions, we know that $p: {\bf R}^m\rightarrow
({\bf R}^n,\omega)$ is continuous if and only if for any number
$\epsilon
>0$, there is  number $\delta=\min\{\delta_1,\delta_2\}$ such that

$$\|\omega(p(\overline{u}))-\omega(p(\overline{v}))\|<\epsilon$$

\no for $\forall\overline{u},\overline{v}\in {\bf R}^m$ with
$\|\overline{u}-\overline{v})\|<\delta$. \ \ $\natural$

\vskip 4mm

\no{\bf Corollary $6.1.1$} \ {\it If $m\geq n+1$, let $\omega:{\bf
R}^{n}\rightarrow{\bf R}^{m-n}$ be a continuous mapping, then
$({\bf R}^n,\omega(p({\bf R}^m)))$ is a pseudo-face of ${\bf R}^m$
in $({\bf R}^n,\omega)$ with}

$$p(x_1,x_2,\cdots,x_n,x_{n+1},\cdots,x_m)=\omega(x_1,x_2,\cdots,x_n).$$

\no{\it Particularly, if $m=3, n=2$ and $\omega$ is an angle
function, then $({\bf R}^n,\omega(p({\bf R}^m)))$ is a pseudo-face
with $p(x_1,x_2,x_3)=\omega(x_1,x_2)$.}

\vskip 3mm

There is a simple relation for a continuous mapping between Euclid
spaces and that of between pseudo-faces established in the next
result.

\vskip 4mm

\no{\bf Theorem $6.1.2$} \ {\it Let $g: {\bf R}^m\rightarrow{\bf
R}^m$ and $p:{\bf R}^m\rightarrow({\bf R}^n,\omega)$ be continuous
mappings. Then $pgp^{-1}: ({\bf R}^n,\omega)\rightarrow({\bf
R}^n,\omega)$ is also a continuous mapping.}

\vskip 3mm

{\it Proof} \ Because a composition of continuous mappings is a
continuous mapping, we know that $pgp^{-1}$ is continuous.

Now for $\forall \omega(x_1,x_2,\cdots,x_n)\in({\bf R}^n,\omega)$,
assume that $p(y_1,y_2,\cdots,y_m)=\omega(x_1,x_2,$ $\cdots,x_n)$,
$g(y_1,y_2,\cdots,y_m)=(z_1,z_2,\cdots,z_m)$ and
$p(z_1,z_2,\cdots,z_m)=\omega(t_1,t_2,\cdots,t_n)$. Then
calculation shows that

\begin{eqnarray*}
pgp^{-}(\omega(x_1,x_2,\cdots,x_n))&=& pg(y_1,y_2,\cdots,y_m)\\
&=& p(z_1,z_2,\cdots,z_m)=\omega(t_1,t_2,\cdots,t_n)\in({\bf
R}^n,\omega).
\end{eqnarray*}

\no Whence, $pgp^{-1}$ is a continuous mapping and  $pgp^{-1}:
({\bf R}^n,\omega)\rightarrow({\bf R}^n,\omega). \ \ \natural$

\vskip 4mm

\no{\bf Corollary $6.1.2$} \ {\it Let $C({\bf R}^m)$ and $C({\bf
R}^n,\omega)$ be sets of continuous mapping on an Euclid space
${\bf R}^m$ and an pseudo-metric space $({\bf R}^n,\omega)$. If
there is a pseudo-space for ${\bf R}^m$ in $({\bf R}^n,\omega)$.
Then there is a bijection between $C({\bf R}^m)$ and $C({\bf
R}^n,\omega)$.}

\vskip 3mm

For a body ${\mathcal B}$ in an Euclid space ${\bf R}^m$, its
shape in a pseudo-face $({\bf R}^n,\omega(p({\bf R}^m)))$ of ${\bf
R}^m$ in $({\bf R}^n,\omega)$ is called a {\it pseudo-shape} of
${\mathcal B}$. We get results for pseudo-shapes of a ball in the
following.

\vskip 4mm

\no{\bf Theorem $6.1.3$} \ {\it Let ${\mathcal B}$ be an
$(n+1)$-ball of radius $R$ in a space ${\bf R}^{n+1}$, i.e.,}

$$x_1^2+x_2^2+\cdots+x_n^2+t^2\leq R^2.$$

\no{\it Define a continuous mapping $\omega:{\bf
R}^n\rightarrow{\bf R}^n$ by}

$$\omega(x_1,x_2,\cdots,x_n)=\varsigma t(x_1,x_2,\cdots,x_n)$$

\no{\it for a real number $\varsigma$ and a continuous mapping
$p:{\bf R}^{n+1}\rightarrow{\bf R}^n$ by}

$$p(x_1,x_2,\cdots,x_n,t)=\omega(x_1,x_2,\cdots,x_n).$$

\no{\it Then the pseudo-shape of ${\mathcal B}$ in $({\bf
R}^n,\omega)$ is a ball of radius $\frac{\sqrt{R^2-t^2}}{\varsigma
t}$ for any parameter $t, -R\leq t\leq R$. Particularly, for the
case of $n=2$ and $\varsigma=\frac{1}{2}$, it is a circle of
radius $\sqrt{R^2-t^2}$ for parameter $t$ and an elliptic ball in
${\bf R}^3$ as shown in Fig.$6.1$.}

\includegraphics[bb=-40 10 100 170]{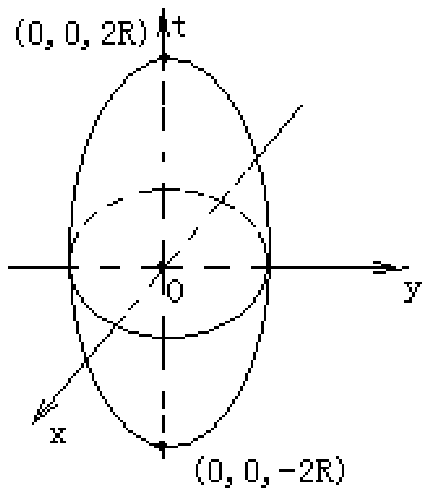}

\vskip 3mm

\c{\bf Fig.$6.1$}

\vskip 3mm

{\it Proof} \ For any parameter $t$, an $(n+1)$-ball

$$x_1^2+x_2^2+\cdots+x_n^2+t^2\leq R^2$$

\no can be transferred to an $n$-ball

$$x_1^2+x_2^2+\cdots+x_n^2\leq R^2-t^2$$

\no of radius $\sqrt{R^2-t^2}$. Whence, if we define a continuous
mapping on ${\bf R}^n$ by

$$\omega(x_1,x_2,\cdots,x_n)=\varsigma t(x_1,x_2,\cdots,x_n)$$

\no and

$$p(x_1,x_2,\cdots,x_n,t)=\omega(x_1,x_2,\cdots,x_n),$$

\no then we get an $n$-ball

$$x_1^2+x_2^2+\cdots+x_n^2\leq \frac{R^2-t^2}{\varsigma^2t^2},$$

\no of ${\mathcal B}$ under $p$ for any parameter $t$, which is
the pseudo-face of ${\mathcal B}$ for a parameter $t$ by
definition.

For the case of $n=2$ and $\varsigma=\frac{1}{2}$, since its
pseudo-face is a circle in an Euclid plane and $-R\leq t\leq R$,
we get an elliptic ball as shown in Fig.$6.1$.\quad\quad
$\natural$

Similarly, if we define $\omega(x_1,x_2,\cdots,x_n)=
2\angle(\overrightarrow{OP},Ot)$ for a point $P=(x_1,x_2,\cdots,$
$x_n,t)$, i.e., an angle function, then we can also get a result
like Theorem $6.1.2$ for these pseudo-shapes of an $(n+1)$-ball.

\vskip 4mm

\no{\bf Theorem $6.1.4$} \ {\it Let ${\mathcal B}$ be an
$(n+1)$-ball of radius $R$ in a space ${\bf R}^{n+1}$, i.e.,}

$$x_1^2+x_2^2+\cdots+x_n^2+t^2\leq R^2.$$

\no{\it Define a continuous mapping $\omega:{\bf
R}^n\rightarrow{\bf R}^n$ by}

$$\omega(x_1,x_2,\cdots,x_n)=2\angle(\overrightarrow{OP},Ot)$$

\no{\it for a point $P$ on ${\mathcal B}$ and a continuous mapping
$p:{\bf R}^{n+1}\rightarrow{\bf R}^n$ by}

$$p(x_1,x_2,\cdots,x_n,t)=\omega(x_1,x_2,\cdots,x_n).$$

\no{\it Then the pseudo-shape of ${\mathcal B}$ in $({\bf
R}^n,\omega)$ is a ball of radius $\sqrt{R^2-t^2}$ for any
parameter $t, -R\leq t\leq R$. Particularly, for the case of
$n=2$, it is a circle of radius $\sqrt{R^2-t^2}$ for parameter $t$
and a body in ${\bf R}^3$ with equations}

$$\oint\arctan(\frac{t}{x})=2\pi \ \ and \ \ \oint\arctan(\frac{t}{y})=2\pi$$

\no{\it for curves of its intersection with planes $XOT$ and
$YOT$.}

\vskip 3mm

{\it Proof} \ The proof is similar to the proof of Theorem
$6.1.3$. For these equations

$$\oint\arctan(\frac{t}{x})=2\pi \ \ {\rm or} \ \ \oint\arctan(\frac{t}{y})=2\pi$$

\no of curves on planes $XOT$ or $YOT$ in the case of $n=2$, they
are implied by the geometrical meaning of an angle function.
\quad\quad $\natural$

For an Euclid space ${\bf R}^n$, we can get a subspace sequence

$${\bf R}_0\supset{\bf R}_1\supset\cdots\supset{\bf R}_{n-1}\supset{\bf R}_n,$$

\no where the dimensional of ${\bf R}_i$ is $n-i$ for $1\leq i\leq
n$  and ${\bf R}_{n}$ is just a point. But we can not get a
sequence reversing the order, i.e., a sequence

$${\bf R}_0\subset{\bf R}_1\subset\cdots\subset{\bf R}_{n-1}\subset{\bf R}_{n}$$

\no in classical space theory. By applying Smarandache
multi-spaces, we can really find this kind of sequence by the next
result, which can be used to explain a well-known model for our
cosmos in M-theory.

\vskip 4mm

\no{\bf Theorem $6.1.5$} \ {\it Let $P=(x_1,x_2,\cdots,x_n)$ be a
point of ${\bf R}^n$. Then there are subspaces of dimensional $s$
in $P$ for any integer $s, 1\leq s\leq n$.}

\vskip  3mm

{\it Proof} \ Notice that in an Euclid space ${\bf R}^n$, there is
a basis $e_1=(1,0,0,\cdots,0)$, $e_2=(0,1,0,\cdots,0)$, $\cdots$,
$e_i=(0,\cdots, 0,1,0,\cdots,0)$ (every entry is $0$ unless the
$i$-th entry is $1$), $\cdots$, $e_n=(0,0,\cdots,0,1)$ such that

$$(x_1,x_2,\cdots,x_n)=x_1e_1+x_2e_2+\cdots+x_ne_n$$

\no for any point $(x_1,x_2,\cdots,x_n)$ of ${\bf R}^n$. Now we
consider a linear space ${\bf R}^- =(V,+_{new},\circ_{new})$ on a
field $F=\{a_i,b_i,c_i,\cdots,$ $d_i; i\geq 1\}$, where

$$V =\{x_1,x_2,\cdots,x_n\}.$$

\no Not loss of generality, we assume that $x_1,x_2,\cdots, x_s$
are independent, i.e., if there exist scalars $a_1,a_2,\cdots,a_s$
such that

$$a_1\circ_{new}x_1+_{new}a_2\circ_{new}x_2+_{new}\cdots+_{new}a_s\circ_{new}x_s=0,$$

\no then $a_1=a_2=\cdots=0_{new}$ and there are scalars $b_i,
c_i,\cdots, d_i$ with $1\leq i\leq s$ in ${\bf R}^-$ such that

$$x_{s+1}=b_1\circ_{new}x_1+_{new}b_2\circ_{new}x_2+_{new}\cdots+_{new}b_s\circ_{new}x_s;$$

$$x_{s+2}=c_1\circ_{new}x_1+_{new}c_2\circ_{new}x_2+_{new}\cdots+_{new}c_s\circ_{new}x_s;$$

$$\cdots\cdots\cdots\cdots\cdots\cdots\cdots\cdots\cdots\cdots;$$

$$x_n=d_1\circ_{new}x_1+_{new}d_2\circ_{new}x_2+_{new}\cdots+_{new}d_s\circ_{new}x_s.$$

\no Therefore, we get a subspace of dimensional $s$ in a point $P$
of ${\bf R}^n$.\quad\quad $\natural$

\vskip 4mm

\no{\bf Corollary $6.1.3$} \ {\it Let $P$ be a point of an Euclid
space ${\bf R}^n$. Then there is a subspace sequence}

$${\bf R}_0^-\subset{\bf R}_1^-\subset\cdots\subset{\bf R}_{n-1}^-\subset{\bf R}_n^-$$

\no{\it such that ${\bf R}_n^-= \{ P\}$ and the dimensional of the
subspace ${\bf R}_i^-$ is $n-i$, where $1\leq i\leq n$.}

\vskip 3mm

{\it Proof} \ Applying Theorem $6.1.5$ repeatedly, we get the
desired sequence. \ \ $\natural$

\vskip 8mm

\no{\bf \S $6.2.$ Relativity Theory}

\vskip 4mm

\no In theoretical physics, these spacetimes are used to describe
various states of particles dependent on the time parameter in an
Euclid space ${\bf R}^3$. There are two kinds of spacetimes. An
{\it absolute spacetime} is an Euclid space ${\bf R}^3$ with an
independent time, denoted by $(x_1,x_2,x_3|t)$ and a {\it relative
spacetime} is an Euclid space ${\bf R}^4$, where time is the
$t$-axis, seeing also in $[30]-[31]$ for details.

A point in a spacetime is called an {\it event}, i.e., represented
by

$$(x_1,x_2,x_3)\in{\bf R}^3 \ \ {\rm and} \ \ t\in{\bf R}^+$$

\no in an absolute spacetime in the Newton's mechanics and

$$(x_1,x_2,x_3,t)\in {\bf R}^4$$

\no with time parameter $t$ in a relative space-time used in the
Einstein's relativity theory.

For two events $A_1=(x_1,x_2,x_3|t_1)$ and
$A_2=(y_1,y_2,y_3|t_2)$, the {\it time interval} $\triangle t$ is
defined by $\triangle t=t_1-t_2$ and the {\it space interval}
$\triangle(A_1,A_2)$ by

$$\triangle(A_1,A_2)=\sqrt{(x_1-y_1)^2+(x_2-y_2)^2+(x_3-y_3)^2}.$$

Similarly, for two events $B_1=(x_1,x_2,x_3,t_1)$ and
$B_2=(y_1,y_2,y_3,t_2)$, the {\it spacetime interval} $\triangle
s$ is defined by

$$\triangle^2 s= -c^2\triangle t^2+\triangle^2(B_1,B_2),$$

\no where $c$ is the speed of the light in vacuum. In Fig.$6.2$, a
spacetime only with two parameters $x,y$ and the time parameter
$t$ is shown.

\includegraphics[bb=10 10 100 140]{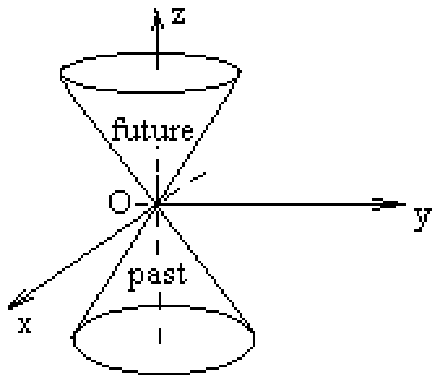}

\vskip 3mm

\c{\bf Fig.$6.2$}

\vskip 2mm

The Einstein's spacetime is an uniform linear space. By the
assumption of linearity of a spacetime and invariance of the light
speed, it can be shown that the invariance of the space-time
intervals, i.e.,\vskip 3mm

{\it For two reference systems $S_1$ and $S_2$ with a homogenous
relative velocity, there must be}

$$\triangle s^2 = \triangle s'^2.$$\vskip 2mm

We can also get the Lorentz transformation of spacetimes or
velocities by this assumption. For two parallel reference systems
$S_1$ and $S_2$, if the velocity of $S_2$ relative to $S_1$ is $v$
along $x$-axis such as shown in Fig.$6.3$,

\includegraphics[bb=10 10 100 120]{sgm92.eps}

\c{\bf Fig.$6.3$}

\vskip 2mm

\no then we know the {\it Lorentz transformation of spacetimes}

\[\left\{
\begin{array}{l}
x_2=\frac{x_1-vt_1}{\sqrt{1-(\frac{v}{c})^2}}\\
y_2=y_1\\
z_2=z_1 \\
t_2=\frac{t_1-\frac{v}{c}x_1}{\sqrt{1-(\frac{v}{c})^2}}
\end{array}
\right.
\]

\no and the {\it transformation of velocities}

\[\left\{
\begin{array}{l}
v_{x_2}=\frac{v_{x_1}-v}{1-\frac{vv_{x_1}}{c^2}}\\
v_{y_2}=\frac{v_{y_1}\sqrt{1-(\frac{v}{c})^2}}{1-\frac{vv_{x_1}}{c^2}}\\
v_{z_2}=
\frac{v_{z_1}\sqrt{1-(\frac{v}{c})^2}}{1-\frac{vv_{x_1}}{c^2}}.
\end{array}
\right.
\]

\vskip 4mm

In the relative spacetime, the {\it general interval} is defined
by

$$ds^2=g_{\mu\nu}dx^{\mu}dx^{\nu},$$

\no where $g_{\mu\nu}=g_{\mu\nu}(x^{\sigma},t)$ is a metric
dependent on the space and time. We can also introduce the
invariance of general intervals, i.e.,

$$ds^2=g_{\mu\nu}dx^{\mu}dx^{\nu}
=g'_{\mu\nu}dx'^{\mu}dx'^{\nu}.$$

\no Then the {\it Einstein's equivalence principle} says that
\vskip 3mm

{\it There are no difference for physical effects of the inertial
force and the gravitation in a field small enough.}\vskip 2mm

An immediately consequence of the Einstein's equivalence principle
is the idea of the {\it geometrization of gravitation}, i.e.,
considering the curvature at each point in a spacetime to be all
effect of gravitation([$18$]), which is called a {\it
gravitational factor} at this point.

Combining these discussions in Section $6.1$ with the Einstein's
idea of the geometrization of gravitation, we get a result for
spacetimes in the theoretical physics.

\vskip 4mm

\no{\bf Theorem $6.2.1$} \ {\it Every spacetime is a pseudo-face
in an Euclid pseudo-space, especially, the Einstein's space-time
is ${\bf R}^n$ in $({\bf R}^4,\omega)$ for an integer $n, n\geq
4$.}

\vskip 3mm

By the uniformity of a spacetime, we get an equation by
equilibrium of vectors in a cosmos.

\vskip 4mm

\no{\bf Theorem $6.2.2$} \ {\it By the assumption of uniformity
for a spacetime in $({\bf R}^4,\omega)$, there exists an
anti-vector $\omega_O^-$ of $\omega_O$ along any orientation
$\overrightarrow{O}$ in ${\bf R}^4$ such that}

$$\omega_O + \omega_O^-=0.$$

\vskip 3mm

{\it Proof} \ Since ${\bf R}^4$ is uniformity, By the principle of
equilibrium in a uniform space, along any orientation
$\overrightarrow{O}$ in ${\bf R}^4$, there must exists an
anti-vector $\omega_O^-$ of $\omega_O$ such that

$$\omega_O + \omega_O^-=0. \ \ \natural$$

Theorem $6.2.2$ has many useful applications. For example, let

$$\omega_{\mu\nu} = R_{\mu\nu}-\frac{1}{2}Rg_{\mu\nu}+\lambda g_{\mu\nu},$$

\no then we know that

$$\omega_{\mu\nu}^-= -8\pi GT_{\mu\nu}.$$

\no in a gravitational field. Whence, we get the {\it Einstein's
equation of gravitational field}

$$R_{\mu\nu}-\frac{1}{2}Rg_{\mu\nu}+\lambda g_{\mu\nu}=-8\pi GT_{\mu\nu}$$

\no by the equation in Theorem $6.2.2$ which is widely used for
our cosmos by physicists. In fact, there are two assumptions for
our cosmos in the following. One is partially adopted from the
Einstein's, another is just suggested by ours.

\vskip 3mm

\no{\bf Postulate $6.2.1$} \ {\it At the beginning our cosmos is
homogenous.}

\vskip 2mm

\no{\bf Postulate $6.2.2$} \ {\it Human beings can only survey
pseudo-faces of our cosmos by observations and experiments.}

\vskip 2mm

Applying these postulates, the Einstein's equation of
gravitational field and the {\it cosmological principle}, i.e.,
{\it there are no difference at different points and different
orientations at a point of a cosmos on the metric $10^4l.y.$}, we
can get a standard model of cosmos, also called the {\it Friedmann
cosmos}, seeing $[18]$,$[26],[28]$,$[30]-[31]$,$[79]$ and $[95]$
for details. In this model, its line element $ds$ is

$$ds^2 = -c^2dt^2+a^2(t)[\frac{dr^2}{1-Kr^2}+r^2(d\theta^2+\sin^2\theta d\varphi^2)]$$

\no and cosmoses are classified into three types:

\vskip 3mm

{\bf Static Cosmos}: \ \ $da/dt \ =0$;

\vskip 2mm

{\bf Contracting Cosmos}: \ \ $da/dt \ <0$;

\vskip 2mm

{\bf Expanding Cosmos}: \ \ $da/dt \ >0$.\vskip 3mm

By the Einstein's view, our living cosmos is the static cosmos.
That is why he added a cosmological constant $\lambda$ in his
equation of gravitational field. But unfortunately, our cosmos is
an expanding cosmos found by Hubble in 1929. As a by-product, the
shape of our cosmos described by S.Hawking in $[30]-[32]$ is
coincide with these results gotten in Section $6.1$.

\vskip 8mm

\no{\bf \S $6.3$ \ A Multi-Space Model for Cosmoses}

\vskip 5mm

\no{\bf $6.3.1.$ What is M-theory}

\vskip 3mm

\no Today, we have know that all matter are made of atoms and
sub-atomic particles, held together by four fundamental forces:
{\it gravity, electro-magnetism, strong nuclear force} and {\it
weak force}. Their features are partially explained by the {\it
quantum theory} and the {\it relativity theory}. The former is a
theory for the microcosm but the later is for the macrocosm.
However, these two theories do not resemble each other in any way.
The quantum theory reduces forces to the exchange of discrete
packet of quanta, while the relativity theory explains the cosmic
forces by postulating the smooth deformation of the fabric
spacetime.

As we known, there are two string theories : the {$E_8\times E_8$
heterotic string}, the {\it SO(32) heterotic string} and three
superstring theories: the {\it SO(32) Type I string}, the {\it
Type $IIA$} and {\it Type $IIB$} in superstring theories. Two
physical theories are {\it dual} to each other if they have
identical physics after a certain mathematical transformation.
There are {\it T-duality} and {\it S-duality} in superstring
theories defined in the following table $6.1$([$15$]).

\vskip 3mm

\begin{center}
\begin{tabular}{|c|c|c|} \hline
$ \ \ $ & fundamental \ string &  dual \ string  \\ \hline
$T-duality$ & $Radius\leftrightarrow 1/(radius)$ & $charge\leftrightarrow 1/(charge)$ \\
$\ \ $ & $Kaluza-Klein\leftrightarrow Winding$ &
$Electric\leftrightarrow Magnet$ \\ \hline
$S-duality$ & $charge\leftrightarrow 1/(charge)$ & $Radius\leftrightarrow 1/(Radius)$  \\
$ \ \ $ & $Electric\leftrightarrow Magnetic$ &
$Kaluza-Klein\leftrightarrow Winding$\\ \hline
\end{tabular}
\end{center}
\vskip 2mm

\c{\bf table $6.1$}

\vskip 2mm

We already know some profound properties for these spring or
superspring theories, such as:

\vskip 3mm

($i$) {\it Type $IIA$ and $IIB$ are related by T-duality, as are
the two heterotic theories}.

($ii$) {\it Type $I$ and heterotic $SO(32)$ are related by
S-duality and Type $IIB$ is also S-dual with itself}.

($iii$) {\it Type $II$ theories have two supersymmetries in the
$10$-dimensional sense, but the rest just one}.

($iv$) {\it Type $I$ theory is special in that it is based on
unoriented open and closed strings, but the other four are based
on oriented closed strings}.

($v$) {\it The $IIA$ theory is special because it is
non-chiral(parity conserving), but the other four are
chiral(parity violating)}.

($vi$) {\it In each of these cases there is an 11th dimension that
becomes large at strong coupling. For substance, in the $IIA$ case
the 11th dimension is a circle and in $IIB$ case it is a line
interval, which makes 11-dimensional spacetime display two
10-dimensional boundaries}.

($vii$) {\it The strong coupling limit of either theory produces
an 11-dimensional spacetime}.

($viii$) $\cdots$, etc..

\vskip 2mm

The M-theory was established by Witten in 1995 for the unity of
those two string theories and three superstring theories, which
postulates that all matter and energy can be reduced to {\it
branes} of energy vibrating in an $11$ dimensional space. This
theory gives us a compelling explanation of the origin of our
cosmos and combines all of existed string theories by showing
those are just special cases of M-theory such as shown in table
$6.2$.

\[ M-theory \left\{
\begin{array}{l}
E_8\times E_8 \ heterotic \ string\\
SO(32)
\ heterotic \ string\\
SO(32) \ Type \ I \ string\\
Type \ IIA\\
Type \ IIB.
\end{array}
\right.
\]

\vskip 3mm \c{\bf Table $6.2$}\vskip 2mm

See Fig.$6.4$ for the M-theory planet in which we can find a
relation of M-theory with these two strings or three superstring
theories. \vskip 3mm

\includegraphics[bb=10 10 240 400]{sgm93.eps}

\vskip 5mm

\c{\bf Fig.$6.4$}

\vskip 2mm

As it is widely accepted that our cosmos is in accelerating
expansion, i.e., our cosmos is most possible an accelerating
cosmos of expansion, it should satisfies the following condition

$$\frac{d^2a}{dt^2} \ > \ 0.$$

The {\it Kasner} type metric

$$ds^2 = -dt^2+a(t)^2d_{{\bf R}^3}^2+b(t)^2ds^2(T^m)$$

\no solves the $4+m$ dimensional vacuum Einstein equations if

$$a(t)=t^{\mu} \ \ {\rm and} \ \ b(t)=t{\nu}$$

\no with

$$\mu=\frac{3\pm\sqrt{3m(m+2)}}{3(m+3)}, \nu=\frac{3\mp\sqrt{3m(m+2)}}{3(m+3)}.$$

These solutions in general do not give an accelerating expansion
of spacetime of dimension $4$. However, by using the time-shift
symmetry

$$t\rightarrow t_{+\infty}-t, \ \ a(t)=(t_{+\infty}-t)^{\mu},$$

\no we see that yields a really accelerating expansion since

$$
\frac{da(t)}{dt} >0 \ {\rm and} \ \frac{d^2a(t)}{dt^2} >0.
$$

According to M-theory, our cosmos started as a perfect $11$
dimensional space with nothing in it. However, this $11$
dimensional space was unstable. The original $11$ dimensional
spacetime finally cracked into two pieces, a $4$ and a $7$
dimensional cosmos. The cosmos made the $7$ of the $11$ dimensions
curled into a tiny ball, allowing the remaining $4$ dimensional
cosmos to inflate at enormous rates. This origin of our cosmos
implies a multi-space result for our cosmos verified by Theorem
$6.1.5$.

\vskip 4mm

\no{\bf Theorem $6.3.1$} \ {\it The spacetime of M-theory is a
multi-space with a warping ${\bf R}^7$ at each point of ${\bf
R}^4$.}

\vskip 3mm

Applying Theorem $6.3.1$, an example for an accelerating expansion
cosmos of $4$-dimensional cosmos from supergravity
compactification on hyperbolic spaces is the {\it
Townsend-Wohlfarth type} in which the solution is

$$ds^2=e^{-m\phi(t)}(-S^6dt^2+S^2dx_3^2)+r_C^2e^{2\phi(t)}ds_{H_m}^2,$$

\no where

$$\phi(t)=\frac{1}{m-1}(\ln K(t)-3\lambda_0t), \ \ \
S^2=K^{\frac{m}{m-1}}e^{-\frac{m+2}{m-1}\lambda_0t}$$

\no and

$$K(t)=\frac{\lambda_0\zeta r_c}{(m-1)\sin[\lambda_0\zeta|t+t_1|]}$$

\no with $\zeta=\sqrt{3+6/m}$. This solution is obtainable from
space-like brane solution and if the proper time $\varsigma$ is
defined by $d\varsigma= S^3(t)dt$, then the conditions for
expansion and acceleration are $\frac{dS}{d\varsigma}>0$ and
$\frac{d^2S}{d\varsigma^2}>0$. For example, the expansion factor
is $3.04$ if $m=7$, i.e., a really expanding cosmos.

\vskip 4mm

\no{\bf $6.3.2.$ A pseudo-face model for $p$-branes}

\vskip 3mm

\no In fact, M-theory contains much more than just strings, which
is also implied in Fig.$6.4$. It contains both higher and lower
dimensional objects, called {\it branes}. A {\it brane} is an
object or subspace which can have various spatial dimensions. For
any integer $p\geq 0$, a {\it $p$-brane} has length in $p$
dimensions, for example, a {\it $0$-brane} is just a point; a {\it
$1$-brane} is a string and a {\it $2$-brane} is a surface or
membrane $\cdots$.

Two branes and their motion have been shown in Fig.$6.5$ where (a)
is a $1$-brane and (b) is a $2$-brane.

\includegraphics[bb=80 10 240 270]{sgm94.eps}

\vskip 2mm

\c{\bf Fig.$6.5$}

\vskip 2mm

Combining these ideas in the pseudo-spaces theory and in M-theory,
a model for ${\bf R}^{m}$ is constructed in the below.

\vskip 4mm

\no{\bf Model $6.3.1$} \ {\it For each $m$-brane ${\bf B}$ of a
space ${\bf R}^{m}$, let $(n_1({\bf B}),n_2({\bf B}),\cdots,
n_p({\bf B}))$ be its unit vibrating normal vector along these $p$
directions and $q:{\bf R}^{m}\rightarrow{\bf R}^4$ a continuous
mapping. Now for $\forall P\in{\bf B}$, define}

$$\omega(q(P))=(n_1(P),n_2(P),\cdots, n_p(P)).$$

\no{\it Then $({\bf R}^4,\omega)$ is a pseudo-face of ${\bf R}^m$,
particularly, if $m=11$, it is a pseudo-face for the M-theory.}

\vskip 3mm

For the case of $p=4$, interesting results are obtained by
applying results in Chapters $5$.

\vskip 4mm

\no{\bf Theorem $6.3.2$} \ {\it For a sphere-like cosmos ${\bf
B}^2$, there is a continuous mapping $q:{\bf B}^2\rightarrow{\bf
R}^2$ such that its spacetime is a pseudo-plane.}

\vskip 3mm

{\it Proof} \ According to the classical geometry, we know that
there is a projection $q:{\bf B}^2\rightarrow{\bf R}^2$ from a
$2$-ball ${\bf B}^2$ to an Euclid plane ${\bf R}^2$, as shown in
Fig.$6.6$.

\includegraphics[bb=0 10 240 160]{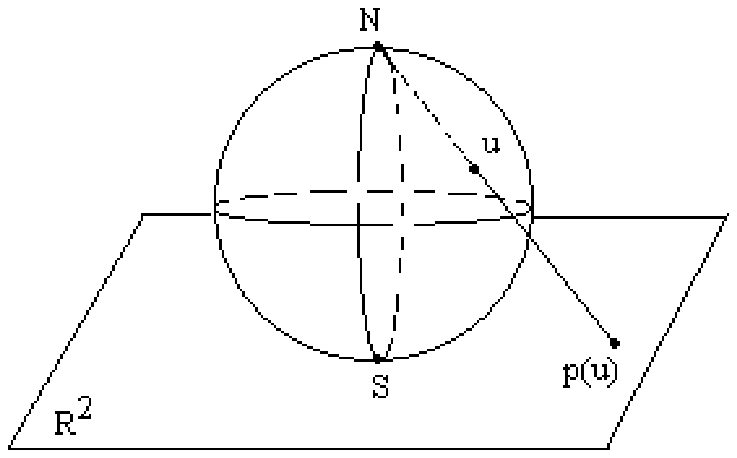}

\vskip 3mm

\c{\bf Fig.$6.6$}

\vskip 2mm

Now for any point $u\in{\bf B}^2$ with an unit vibrating normal
vector $(x(u),y(u),z(u))$, define

$$\omega(q(u))=(z(u),t),$$

\no where $t$ is the time parameter. Then $({\bf R}^2,\omega)$ is
a pseudo-face of $({\bf B}^2,t)$. \ \ $\natural$

Generally, we can also find pseudo-surfaces as a pseudo-face of
sphere-like cosmoses.

\vskip 4mm

\no{\bf Theorem $6.3.3$} \ {\it For a sphere-like cosmos ${\bf
B}^2$ and a surface $S$, there is a continuous mapping $q:{\bf
B}^2\rightarrow S$ such that its spacetime is a pseudo-surface on
$S$.}

\vskip 3mm

{\it Proof} \ According to the classification theorem of surfaces,
an surface $S$ can be combinatorially represented by a
$2n$-polygon for an integer $n, n\geq 1$. If we assume that each
edge of this polygon is at an infinite place, then the projection
in Fig.$6.6$ also enables us to get a continuous mapping $q:{\bf
B}^2\rightarrow S$. Thereby we get a pseudo-face on $S$ for the
cosmos ${\bf B}^2$. \ \ $\natural$

Furthermore, we can construct a combinatorial model for our cosmos
by applying materials in Section $2.5$.

\vskip 4mm

\no{\bf Model $6.3.2$} \ {\it For each $m$-brane ${\bf B}$ of a
space ${\bf R}^{m}$, let $(n_1({\bf B}),n_2({\bf B}),\cdots,
n_p({\bf B}))$ be its unit vibrating normal vector along these $p$
directions and $q:{\bf R}^{m}\rightarrow{\bf R}^4$ a continuous
mapping. Now construct a graph phase $({\mathcal
G},\omega,\Lambda)$ by}

$$V({\mathcal G})=\{p-branes \ q({\bf B})\},$$

$$E({\mathcal G})=\{(q({\bf B}_1),q({\bf B}_2))| there \ is \ an \ action \ between \ {\bf B}_1 \
and \ {\bf B}_2\},$$

$$\omega(q({\bf B}))=(n_1({\bf B}),n_2({\bf B}),\cdots,
n_p({\bf B})),$$

\no {\it and}

$$\Lambda(q({\bf B}_1),q({\bf B}_2))= \ forces \ between \ {\bf B}_1 \ and \ {\bf B}_2.$$

\no{\it Then we get a graph phase $({\mathcal G},\omega,\Lambda)$
in ${\bf R}^4$. Similarly, if $m=11$, it is a graph phase for the
M-theory.}

\vskip 3mm

If there are only finite $p$-branes in our cosmos, then Theorems
$6.3.2$ and $6.3.3$ can be restated as follows.

\vskip 4mm

\no{\bf Theorem $6.3.4$} \ {\it For a sphere-like cosmos ${\bf
B}^2$ with finite $p$-branes and a surface $S$, its spacetime is a
map geometry on $S$.}

\vskip 3mm

Now we consider the transport of a graph phase $({\mathcal
G},\omega,\Lambda)$ in ${\bf R}^m$ by applying results in Sections
$2.3$ and $2.5$.

\vskip 4mm

\no{\bf Theorem $6.3.5$} \ {\it A graph phase $({\mathcal
G}_1,\omega_1,\Lambda_1)$ of space  ${\bf R}^m$ is transformable
to a graph phase $({\mathcal G}_2,\omega_2,\Lambda_2)$ of space
${\bf R}^n$ if and only if ${\mathcal G}_1$ is embeddable in ${\bf
R}^n$ and there is a continuous mapping $\tau$ such that
$\omega_2=\tau(\omega_1)$ and $\Lambda_2=\tau(\Lambda_1)$.}

\vskip 3mm

{\it Proof} \ By the definition of transformations, if $({\mathcal
G}_1,\omega_1,\Lambda_1)$ is transformable to $({\mathcal
G}_2,\omega_2,\Lambda_2)$, then there must be ${\mathcal G}_1$ is
embeddable in ${\bf R}^n$ and  there is  a continuous mapping
$\tau$ such that $\omega_2=\tau(\omega_1)$ and
$\Lambda_2=\tau(\Lambda_1)$.

Now if ${\mathcal G}_1$ is embeddable in ${\bf R}^n$ and there is
a continuous mapping $\tau$ such that $\omega_2=\tau(\omega_1)$,
$\Lambda_2=\tau(\Lambda_1)$, let $\varsigma: {\mathcal
G}_1\rightarrow{\mathcal G}_2$ be a continuous mapping from
${\mathcal G}_1$ to ${\mathcal G}_2$, then $(\varsigma,\tau)$ is
continuous and

$$(\varsigma,\tau): ({\mathcal
G}_1,\omega_1,\Lambda_1)\rightarrow ({\mathcal
G}_2,\omega_2,\Lambda_2).$$

\no Therefore $({\mathcal G}_1,\omega_1,\Lambda_1)$ is
transformable to $({\mathcal G}_2,\omega_2,\Lambda_2)$. \ \
$\natural$

Theorem $6.3.5$ has many interesting consequences as by-products.

\vskip 4mm

\no{\bf Corollary $6.3.1$} \ {\it A graph phase $({\mathcal
G}_1,\omega_1,\Lambda_1)$ in ${\bf R}^m$ is transformable to a
planar graph phase $({\mathcal G}_2,\omega_2,\Lambda_2)$ if and
only if ${\mathcal G}_2$ is a planar embedding of ${\mathcal G}_1$
and there is a continuous mapping $\tau$ such that
$\omega_2=\tau(\omega_1)$, $\Lambda_2=\tau(\Lambda_1)$ and vice
via, a planar graph phase $({\mathcal G}_2,\omega_2,\Lambda_2)$ is
transformable to a graph phase $({\mathcal
G}_1,\omega_1,\Lambda_1)$ in ${\bf R}^m$ if and only if ${\mathcal
G}_1$ is an embedding of ${\mathcal G}_2$ in ${\bf R}^m$ and there
is a continuous mapping $\tau^{-1}$ such that
$\omega_1=\tau^{-1}(\omega_2)$, $\Lambda_1=\tau^{-1}(\Lambda_2)$.}

\vskip 4mm

\no{\bf Corollary $6.3.2$} \ {\it For a continuous mapping $\tau$,
a graph phase $({\mathcal G}_1,\omega_1,\Lambda_1)$ in ${\bf R}^m$
is transformable to a graph phase $({\mathcal
G}_2,\tau(\omega_1),\tau(\Lambda_1))$ in ${\bf R}^n$ with $m,n\geq
3$.}

\vskip 3mm

{\it Proof} \ This result follows immediately from Theorems
$2.3.2$ and $6.3.5$. \ \ $\natural$

This theorem can be also used to explain the problems of {\it
travelling between cosmoses} or {\it getting into the heaven or
hell} for a person. For example, water will go from a liquid phase
to a steam phase by heating and then will go to a liquid phase by
cooling because its phase is transformable between the steam phase
and the liquid phase. For a person on the earth, he can only get
into the heaven or hell after death because the dimension of the
heaven is more than $4$ and that of the hell is less than $4$ and
there does not exist a transformation for an alive person from our
cosmos to the heaven or hell by the biological structure of his
body. Whence, if black holes are really these tunnels between
different cosmoses, the destiny for a cosmonaut unfortunately fell
into a black hole is only the death ($[30][32]$). Perhaps, there
are really other kind of beings in cosmoses or mankind in the
further who can freely change from one phase in a space ${\bf
R}^m$ to another in ${\bf R}^n$ with $m\not=n$, then the
travelling between cosmoses is possible for those beings or
mankind in that time.

\vskip 4mm

\no{\bf $6.3.3.$ A multi-space model of cosmos}

\vskip 4mm

\no Until today, many problems in cosmology are puzzling one's
eyes. Comparing with these vast cosmoses, human beings are very
tiny. In spite of this depressed fact, we can still investigate
cosmoses by our deeply thinking. Motivated by this belief, a
multi-space model for cosmoses is constructed in the following.

\vskip 4mm

\no{\bf Model $6.3.3$} \ {\it A mathematical cosmos is constructed
by a triple $(\Omega,\Delta,T)$, where }

$$\Omega=\bigcup\limits_{i\geq 0}\Omega_i, \  \ \ \Delta=\bigcup\limits_{i\geq 0}O_i$$

\no{\it and $T=\{t_i; i\geq 0\}$ are respectively called the
cosmos, the operation or the time set with the following
conditions hold.}\vskip 2mm

($1$) \ {\it $(\Omega,\Delta)$ is a Smarandache multi-space
dependent on $T$, i.e., the cosmos $(\Omega_i,O_i)$ is dependent
on the time parameter $t_i$ for any integer $i, i\geq 0$.}

($2$) \ {\it For any integer $i, i\geq 0$, there is a sub-cosmos
sequence}

$$(S): \ \Omega_i\supset\cdots\supset\Omega_{i1}\supset\Omega_{i0}$$

\no{\it in the cosmos $(\Omega_i,O_i)$ and for two sub-cosmoses
$(\Omega_{ij},O_i)$ and $(\Omega_{il}, O_i)$, if
$\Omega_{ij}\supset\Omega_{il}$, then there is a homomorphism
$\rho_{\Omega_{ij},\Omega_{il}}:(\Omega_{ij},O_i)\rightarrow(\Omega_{il},O_i)$
such that}\vskip 2mm

($i$) \ {\it for $\forall
(\Omega_{i1},O_i),(\Omega_{i2},O_i)(\Omega_{i3},O_i)\in(S)$, if
$\Omega_{i1}\supset\Omega_{i2}\supset\Omega_{i3}$, then}

$$\rho_{\Omega_{i1},\Omega_{i3}}=
\rho_{\Omega_{i1},\Omega_{i2}}\circ\rho_{\Omega_{i2},\Omega_{i3}},$$

\no{\it where ¡°$\circ$¡± denotes the composition operation on
homomorphisms.}

($ii$) \ {\it for $\forall g,h\in\Omega_i$, if for any integer
$i$, $\rho_{\Omega,\Omega_i}(g)=\rho_{\Omega,\Omega_i}(h)$, then
$g=h$.}

($iii$) \ {\it for $\forall i$, if there is an $f_i\in\Omega_i$
with}

$$\rho_{\Omega_i,\Omega_i\bigcap\Omega_j}(f_i)=
\rho_{\Omega_j,\Omega_i\bigcap\Omega_j}(f_j)$$

\no{\it for integers $i,j,\Omega_i\bigcap\Omega_j\not=\emptyset$,
then there exists an $f\in\Omega$ such that
$\rho_{\Omega,\Omega_i}(f)=f_i$ for any integer $i$.}

\vskip 3mm

Notice that this model is a multi-cosmos model. In the Newton's
mechanics, the Einstein's relativity theory or the M-theory, there
is just one cosmos $\Omega$ and these sub-cosmos sequences are

$${\bf R}^3\supset{\bf R}^2\supset{\bf R}^1\supset{\bf R}^0=\{P\},$$

$${\bf R}^4\supset{\bf R}^3\supset{\bf R}^2\supset{\bf R}^1\supset{\bf R}^0=\{P\}$$

\no and

$${\bf R}^4\supset{\bf R}^3\supset
{\bf R}^2\supset{\bf R}^1\supset{\bf R}^0=\{P\}\supset{\bf
R}_7^-\supset\cdots\supset{\bf R}_1^-\supset{\bf R}_0^-=\{Q\}.$$

These conditions in $(2)$ are used to ensure that a mathematical
cosmos posses a {\it general structure sheaf} of a topological
space, for instance if we equip each multi-space $(\Omega_i,O_i)$
with an abelian group $G_i$ for an integer $i, i\geq 0$, then we
get a structure sheaf on a mathematical cosmos. For {\it general
sheaf theory}, one can see in the reference $[29]$ for details.
This structure enables that a being in a cosmos of higher
dimension can supervises those in lower dimension.

Motivated by this multi-space model of cosmos, we present a number
of conjectures on cosmoses in the following. The first is on the
number of cosmoses and their dimension.

\vskip 4mm

\no{\bf Conjecture $6.3.1$} \ {\it There are infinite many
cosmoses and all dimensions of cosmoses make up an integer
interval $[1,+\infty]$.}

\vskip 3mm

A famous proverbs in Chinese says that {\it seeing is believing
but hearing is unbelieving}, which is also a dogma in the
pragmatism. Today, this view should be abandoned by a
mathematician if he wish to investigate the $21$st mathematics. On
the first, we present a conjecture on the problem of travelling
between cosmoses.

\vskip 4mm

\no{\bf Conjecture $6.3.2$} \ {\it There must exists a kind of
beings who can get from one cosmos into another. There must exists
a kind of being who can goes from a space of higher dimension into
its subspace of lower dimension, especially, on the earth.}

\vskip 3mm

Although nearly every physicist acknowledges the existence of
black holes, those holes are really found by mathematical
calculation. On the opposite, we present the next conjecture.

\vskip 4mm

\no{\bf Conjecture $6.3.3$} \ {\it Contrary to black holes, there
are also white holes at where no matters can arrive including the
light in our cosmos.}

\vskip 3mm

\no{\it Conjecture $6.3.4$} \ {\it Every black hole is also a
white hole in a cosmos.}

\vskip 3mm

Our cosmonauts is good luck if Conjecture $6.3.4$ is true since
they do not need to worry about attracted by these black holes in
our cosmos. Today, a very important task in theoretical and
experimental physics is looking for dark matters. However, we do
not think this would be success by the multi-model of cosmoses.
This is included in the following conjecture.

\vskip 4mm

\no{\bf Conjecture $6.3.5$} \ {\it One can not find dark matters
by experiments since they are in spatial can not be found by human
beings.}

\vskip 3mm

Few consideration is on the relation of the dark energy with dark
matters. But we believe there exists a relation between the dark
energy and dark matters such as stated in the next conjecture.

\vskip 4mm

\no{\bf Conjecture $6.3.6$} \ {\it Dark energy is just the effect
of dark matters.}

\vskip 3mm

\newpage

\thispagestyle{empty}

\thispagestyle{empty} \pagestyle{myheadings} \topmargin 5mm
\headheight 8mm \headsep 10mm \markboth{\scriptsize Linfan Mao: \
Smarandache Multi-Spaces Theory}{\scriptsize References}
{\ \ \ }

\vskip 35mm

\no{\bf References} \vskip 5mm

\re{[1]}A.D.Aleksandrov and V.A.Zalgaller, {\it Intrinsic Geometry
of Surfaces}, American Mathematical Society, Providence, Rhode
Island, 1976.

\re{[2]}I.Antoniadis, Physics with large extra dimensions: String
theory under experimental test, {\it Current Science}, Vol.81,
No.12, 25(2001),1609-1613

\re{[3]}V.L.Arnold, {\it Ordinary differential Equations},
Massachusetts Institute of Technology, 1973.

\re{[4]}P.Bedrossian, G.Chen, R.Schelp, A generalization of Fan's
condition for homiltonicity, pancyclicity and hamiltonian
connectedness, {\it Discrete Math.}, Vol.115, 1993,39-50.

\re{[5]}J.C.Bermond, Hamiltonian graphs, {\it Selected Topics in
Graph Theory }, Academic Press, London (1978).

\re{[5]}N.L.Biggs and A.T.White, {\it Permutation Groups and
Combinatoric Structure}, Cambridge University Press, 1979.

\re{[6]}G.Birkhoff and S.Mac Lane, {\it A Survey of Modern
Algebra}, Macmillan Publishing Co., Inc, 1977.

\re{[7]}B.Bollob$\acute{a}$s, {\it Random Graphs}, Academic Press
Inc. (London) Ltd, 1985.

\re{[8]}J.A.Bondy and V.Chv$\acute{a}$tal, A method in graph
theory, {\it Discrete Math.} 15(1976), 111-136.

\re{[9]}J.A.Bondy and U.S.R.Murty, {\it Graph Theory with
Applications}, The Macmillan Press Ltd, 1976.

\re{[10]}M.Capobianco and J.C.Molluzzo, {\it Examples and
Counterexamples in Graph Theory}, Elsevier North-Holland Inc,
1978.

\re{[11]}G.Chartrand and L.Lesniak, {\it Graphs \& Digraphs},
Wadsworth, Inc., California, 1986.

\re{[12]}W.K.Chen, {\it Applied Graph Theory--Graphs and
Electrical Theory}, North-Holland Publishing Company, New York,
1976.

\re{[13]}S.S.Chern and W.H.Chern, {\it Lectures in Differential
Geometry} (in Chinese), Peking University Press, 2001.

\re{[14]}D.Deleanu, {\it A Dictionary of Smarandache Mathematics},
Buxton University Press, London \& New York,2

\re{[15]}M.J.Duff, A layman's guide to M-theory, {\it arXiv}:
hep-th/9805177, v3, 2 July(1998).

\re{[16]}J.Edmonds, A combinatorial representation for polyhedral
surfaces, {\it Notices Amer. Math. Soc.}, 7(1960).

\re{[17]}G.H.Fan, New sufficient conditions for cycle in graphs,
{\it J.Combinat.Theory}, Ser.B(37),1984,221-227.

\re{[18]}B.J.Fei, {\it relativity theory and Non-Euclid
Geometries}, Science Publisher Press, Beijing, 2005.

\re{[19]}C.Frenk, Where is the missing matter?

\ \ \ \ \ http://hom.flash.net/~csmith0/missing.htm.

\re{[20]}M.Gleiser, Where does matter come from?

\ \ \ \ \ http://hom.flash.net/~csmith0/matter.htm.

\re{[21]}C.Godsil and G.Royle, {\it Algebraic Graph Theory},
Springer-Verlag New York,Inc., 2001.

\re{[22]}J.E.Graver and M.E.Watkins, {\it Combinatorics with
Emphasis on the Theory of Graphs}, Springer-Verlag, New York
Inc,1977.

\re{[23]}J.L.Gross and T.W.Tucker, {\it Topological Graph Theory},
John Wiley \& Sons, 1987.

\re{[24]}B.Gr\"{u}mbaum, {\it Convex Polytopes}, London
Interscience Publishers, 1967.

\re{[25]} B.Gr\"{u}mbaum, Polytopes, graphs and complexes,{\it
Bull.Amer.Math.Soc}, 76(1970), 1131-1201.

\re{[26]}B.Gr\"{u}mbaum, Polytopal graphs, {\it Studies in
Mathematics (Mathematical Association Of America)} ,vol. 12(1975),
201-224.

\re{[27]}A.Guth, An eternity of bubbles,

\ \ \ \ \ http://hom.flash.net/~csmith0/bubbles.htm.

\re{[28]]}J.Hartle, Theories of everything and Hawking's wave
function of the universe, in Gibbones, Shellard and Rankin
eds:{\it The Future of Theoretical Physics and Cosmology},
Cambridge University Press, 2003.

\re{[29]}R.Hartshorne, {\it Algebraic Geometry}, Springer-Verlag
New York, Inc., 1977.

\re{[30]}S.Hawking, {\it A Brief History of Times}, A Bantam
Books/ November, 1996.

\re{[31]}S.Hawking, {\it The Universe in a Nutshell}, A Bantam
Books/ November, 2001.

\re{[32]}S.Hawking, Sixty years in a nutshell, in Gibbones,
Shellard and Rankin eds:{\it The Future of Theoretical Physics and
Cosmology}, Cambridge University Press, 2003.

\re{[33]}K.Hoffman and R.Kunze, {\it Linear Algebra} (Second
Edition), Prentice-Hall, Inc., Englewood Cliffs, New Jersey, 1971.

\re{[34]}L.G.Hua, {\it Introduction to Number Theory} (in
Chinese), Science Publish Press, Beijing ,1979.

\re{[35]}H.Iseri, {\it Smarandache Manifolds}, American Research
Press, Rehoboth, NM,2002.

\re{[36]}H.Iseri, {\it Partially Paradoxist Smarandache
Geometries}, http://www.gallup.unm.
edu/\~smarandache/Howard-Iseri-paper.htm.

\re{[77]}M.Kaku, {\it Hyperspace: A Scientific Odyssey through
Parallel Universe, Time Warps and 10th Dimension}, Oxford Univ.
Press.

\re{[38]}L.Kuciuk and M.Antholy, An Introduction to Smarandache
Geometries, {\it Mathematics Magazine, Aurora, Canada},
Vol.12(2003), and online:

\c{http://www.mathematicsmagazine.com/1-2004/Sm\_Geom\_1\_2004.htm;}

also at New Zealand Mathematics Colloquium, Massey University,
Palmerston

North, New Zealand, December 3-6,2001

\c{http://atlas-conferences.com/c/a/h/f/09.htm;}

also at the International Congress of Mathematicians (ICM2002),
Beijing, China,

20-28, August, 2002,

\c{http://www.icm2002.org.cn/B/Schedule.Section04.htm.}

\re{[39]}J.M.Lee, {\it Riemannian Manifolds}, Springer-Verlag New
York, Inc., 1997.

\re{[40]}E.V.Linder, Seeing darkness: the new cosmology, {\it
arXiv}: astro-ph/0511197 v1, (7) Nov. 2005.

\re{[41]}S.D.Liu et.al, {\it Chaos and Fractals in Natural
Sciences} (in Chinese), Peking University Press, 2003.

\re{[42]}Y.P.Liu, {\it Embeddability in Graphs}, Kluwer Academic
Publishers, Dordrecht/ Boston/ London, 1995.

\re{[43]}Y.P.Liu, {\it Enumerative Theory of Maps}, Kluwer
Academic Publishers, Dordrecht/ Boston/ London, 1999.

\re{[44]}Y.P.Liu, {\it Advances in Combinatorial Maps} (in
Chinese), Northern Jiaotong University Publisher, Beijing, 2003.

\re{[45]}A.Malnic, Group action,coverings and lifts of
automorphisms, {\it Discrete Math}, 182(1998),203-218.

\re{[46]}A.Malinic,R.Nedela and M.$\check{S}$koviera, Lifting
graph automorphisms by voltage assignment, {\it Europ.
J.Combinatorics}, 21(2000),927-947.

\re{[47]}L.F.Mao, The maximum size of self-centered graphs with a
given radius (in Chinese), {J. Xidan University}, Vol. 23(supp),
6-10, 1996.

\re{[48]}L.F.Mao, Hamiltonian graphs with constraints on the
vertices degree in a subgraphs pair, {\it J. Institute of Taiyuan
Machinery}, Vol.15(Supp.), 1994,79-90.

\re{[49]}L.F.Mao, A localization of Dirac's theorem for
hamiltonian graphs, {\it J.Math.Res. \& Exp.}, Vol.18, 2(1998),
188-190.

\re{[50]}L.F.Mao, An extension of localized Fan's condition (in
Chinese), {\it J.Qufu Normal University}, Vol.26, 3(2000), 25-28.

\re{[51]}L.F.Mao, On the panfactorial property of Cayley graphs,
{\it J.Math.Res. $\&$ Exp.}, Vol 22,No.3(2002),383-390.

\re{[52]}L.F.Mao, Localized neighborhood unions condition for
hamiltonian graphs, {\it J.Henan Normal University}(Natural
Science), Vol.30, 1(2002), 16-22.

\re{[53]}L.F.Mao, {\it A census of maps on surface with given
underlying graphs}, A doctor thesis in Northern Jiaotong
University, Beijing, 2002.

\re{[54]}L.F.Mao, Parallel bundles in planar map geometries,
{Scientia Magna}, Vol.1(2005), No.2, 120-133.

\re{[55]}L.F.Mao, {\it On Automorphisms groups of Maps, Surfaces
and Smarandache geometries}, A post-doctoral report at the Chinese
Academy of Sciences(2005), also appearing in {\it Sientia Magna},
Vol.$1$(2005), No.$2$, 55-73.

\re{[56]}L.F.Mao, {\it Automorphism Groups of Maps, Surfaces and
Smarandache Geometries}, American Research Press, 2005.

\re{[57]}L.F.Mao, A new view of combinatorial maps by
Smarandache's notion, {\it arXiv}: math.GM/0506232.

\re{[58]}L.F.Mao, On algebraic multi-group spaces, {\it
arXiv}:math.GM/05010427.

\re{[59]}L.F.Mao, On algebraic multi-ring spaces, {\it
arXiv}:math.GM/05010478.

\re{[60]}L.F.Mao, On algebraic multi-vector spaces, {\it
arXiv}:math.GM/05010479.

\re{[61]}L.F.Mao, On multi-metric spaces, {\it
arXiv}:math.GM/05010480.

\re{[62]}L.F.Mao, An introduction to Smarandache geometries on
maps, Reported at the {\it 2005 International Conference on Graphs
and Combinatorics}, Zhejiang, P.R.China, {\it Scientia
Magna}(accepted).

\re{[63]}L.F.Mao and F.Liu, New localized condition with
$d(x,y)\geq 2$ for hamiltonian graphs (in Chinese), {\it J.Henan
Normal University}(Natural Science), Vol.31, 1(2003), 17-21.

\re{[64]}L.F.Mao and Y.P.Liu, On the eccentricity value sequence
of a simple graph, {\it J.Henan Normal University}(Natural
Science), Vol.29, 4(2001),13-18.

\re{[65]}L.F.Mao and Y.P.Liu, New localized condition for
hamiltonian graphs (in Chinese),  {\it J.Qufu Normal University},
Vol. 27, 2(2001), 18-22.

\re{[66]}L.F.Mao and Y.P.Liu, Group action approach for
enumerating maps on surfaces,{\it J.Applied Math. \& Computing},
vol.13, No.1-2,201-215.

\re{[67]} L.F.Mao, Y.P.Liu, E.L.Wei, The semi-arc automorphism
group of a graph with application to map enumeration, {\it Graphs
and Combinatorics}(in pressing).

\re{[68]}W.S.Massey, {\it Algebraic topology: an introduction},
Springer-Verlag,New York, etc., 1977.

\re{[69]}B.Mohar and C.Thomassen, {\it Graphs on Surfaces}, The
Johns Hopkins University Press, London, 2001.

\re{[70]}G.W.Moore, What is a brane? {\it Notices of the AMS},
Vol.52, No.2(2005), 214-215.

\re{[71]}R.Nedela and M $\check{S}$koviera, Regular embeddings of
canonical double coverings of graphs, {\it J.combinatorial
Theory},Ser.B 67, 249-277(1996).

\re{[72]}R.Nedela and M.$\check{S}$koviera, Regular maps from
voltage assignments and exponent groups, {\it
Europ.J.Combinatorics}, 18(1997),807-823.

\re{[73]}L.Z.Nie and S.S.Ding, {\it Introduction to Algebra} (in
Chinese), Higher Education Publishing Press, 1994.

\re{[74]}J.A.Nieto, Matroids and p-branes, {\it Adv. Theor. Math.
Phys.}, 8(2004), 177-188.

\re{[75]}V.V.Nikulin and I.R.Shafarevlch, {\it Geometries and
Groups}, Springer-Verlag Berlin Heidelberg, 1987.

\re{[76]}R.Penrose, The problem of spacetime singularities
:implication for quantum gravity, in Gibbones, Shellard and Rankin
eds:{\it The Future of Theoretical Physics and Cosmology},
Cambridge University Press, 2003.

\re{[77]}PlanetMath, {\it Smarandache Geometries},
http://planetmath.org/encyclopedia/ SmarandacheGeometries.htm1.

\re{[78]}K.Polthier and M.Schmies, Straightest geodesics on
polyhedral surfaces, in {\it Mathematical Visualization} (ed. by
H.C.Hege and K.Polthier), Springer-Verlag, Berlin, 1998.

\re{[79]}M.Rees, Our complex cosmos and its future, in Gibbones,
Shellard and Rankin eds:{\it The Future of Theoretical Physics and
Cosmology}, Cambridge University Press, 2003.

\re{[80]}D.Reinhard, {\it Graph Theory}, Springer-Verlag New
York,Inc., 2000.

\re{[81]}G.Ringel, {\it Map Color Theorem}, Springer-Verlag,
Berlin, 1974.

\re{[82]}D.J.S.Robinson, {\it A Course in the Theory of Groups},
Springer-Verlag, New York Inc., 1982.

\re{[83]}H.J.Ryser, {\it Combinatorial Mathematics}, The
Mathematical Association of America, 1963.

\re{[84]}R.H.Shi, $2$-neighborhood and hamiltonian conditions,
{\it J.Graph Theory}, Vol.16, 1992,267-271.

\re{[85]}F.Smarandache, {\it A Unifying Field in Logics.
Neutrosopy: Neturosophic Probability, Set, and Logic}, American
research Press, Rehoboth, 1999.

\re{[86]}F.Smarandache, Mixed noneuclidean geometries, {\it eprint
arXiv: math/0010119}, 10/2000.

\re{[87]}F.Smarandache, Neutrosophy, a new Branch of Philosophy,
{\it Multi-Valued Logic}, Vol.8, No.3(2002)(special issue on
Neutrosophy and Neutrosophic Logic), 297-384.

\re{[88]}F.Smarandache, A Unifying Field in Logic: Neutrosophic
Field, {\it Multi-Valued Logic}, Vol.8, No.3(2002)(special issue
on Neutrosophy and Neutrosophic Logic), 385-438.

\re{[89]}F.Smarandache, There is no speed barrier for a wave phase
nor for entangle particles, {\it Progress in Physics}, Vol.1,
April(2005), 85-86.

\re{[90]}L.Smolin, What is the future of cosmology?

\ \ \ \ \ http://hom.flash.net/~csmith0/future.htm.

\re{[91]}S.Stahl, Generalized embedding schemes, {\it J.Graph
Theory}, No.2,41-52, 1978.

\re{[92]}J.Stillwell, {\it Classical topology and combinatorial
group theory}, Springer-Verlag New York Inc., 1980.

\re{[93]}R.Thomas and X.X.Yu, 4-connected projective graphs are
hamiltonian, {\it J.Combin. Theory Ser.B}, 62(1994),114-132.

\re{[94]}C.Thomassen, infinite graphs, in {\it Selected Topics in
Graph Theory}(edited by Lowell W.Beineke and Robin J. Wilson),
Academic Press Inc. (London) LTD, 1983.

\re{[95]}K.Thone, Warping spacetime, in Gibbones, Shellard and
Rankin eds:{\it The Future of Theoretical Physics and Cosmology},
Cambridge University Press, 2003.

\re{[96]}W.Thurston, On the geometry and dynamic of
diffeomorphisms of surfaces, {\it Bull. Amer. Math. Soc.},
19(1988), 417-431.

\re{[97]}W.T.Tutte, A theorem on planar graphs, {\it Trans, Amer.
Math. Soc.}, 82(1956), 99-116.

\re{[98]}W.T.Tutte, What is a maps? in {\it New Directions in the
Theory of Graphs} (ed.by F.Harary), Academic Press, 1973, 309-325.

\re{[99]}W.T.Tutte, Bridges and hamiltonian circuits in planar
graphs. {\it Aequationes Mathematicae}, 15(1977),1-33.

\re{[100]}W.T.Tutte, {\it Graph Theory}, Combridge University
Press, 2001.

\re{[101]]}W.B.Vasantha Kandasamy, {\it Bialgebraic structures and
Smarandache bialgebraic structures}, American Research Press,
2003.

\re{[102]}W.B.Vasantha Kandasamy and F.Smarandache, {\it Basic
Neutrosophic Algebraic Structures and Their Applications to Fuzzy
and Neutrosophic Models}, Hexis, Church Rock, 2004.

\re{[103]}W.B.Vasantha Kandasamy and F.Smarandache, {\it
N-Algebraic Structures and S-N-Algebraic Structures}, HEXIS,
Phoenix, Arizona, 2005.

\re{[104]}H.Walther, {\it Ten Applications of Graph Theory},
D.Reidel Publishing Company, Dordrecht/Boston/Lancaster, 1984.

\re{[105]}C.V.Westenholz, {\it Differential Forms in Mathematical
Physics}, North-Holland Publishing Company. New York, 1981.

\re{[106]}A.T.White, {\it Graphs of Group on Surfaces-
interactions and models}, Elsevier Science B.V. 2001.

\re{[107]}Q.Y.Xong, {\it Ring Theory} (in Chinese), Wuhan
University Press, 1993.

\re{[108]}M.Y.Xu, {\it Introduction to Group Theory} (in
Chinese)(I)(II), Science Publish Press, Beijing ,1999.

\re{[109]}Y.C.Yang and L.F.Mao, On the cyclic structure of
self-centered graphs with $r(G)=3$ (in Chinese), {\it Pure and
Applied Mathematics}, Vol.10(Special Issue), 88-98, 1994.

\re{[110]}H.P.Yap, {\it Some Topics in Graph Theory}, Cambridge
University Press, 1986.

\end{document}